\documentclass{birkjour}
\usepackage[all]{xy}
\usepackage{wrapfig}

\usepackage{epic}
\usepackage{amsmath}[1996/11/01]
\usepackage{amssymb,amsthm,amsfonts,latexsym,epsfig,graphics}
 \def\beql#1#2\eeql{\begin{equation}\label{#1}#2\end{equation}}

\DeclareMathOperator{\Chi}{{\mathcal X}}
\DeclareMathOperator{\Rho}{{\mathcal R}}

\DeclareMathOperator{\Irr}{Irr}

\DeclareMathOperator{\End}{End}
\DeclareMathOperator{\Aut}{Aut}

\DeclareMathOperator{\GL}{GL}

\DeclareMathOperator{\diag}{diag}

\DeclareMathOperator{\id}{id}
\DeclareMathOperator{\SL}{SL}
\DeclareMathOperator{\PSL}{PSL}

\theoremstyle{plain}
\newtheorem{theorem}{Theorem}
\newtheorem{lemma}[theorem]{Lemma}

\newtheorem{proposition}[theorem]{Proposition}
\newtheorem{corollary}[theorem]{Corollary}
\newtheorem{definition}[theorem]{Definition}
\theoremstyle{remark}

\newtheorem{remark}[theorem]{Remark}
\newtheorem{example}[theorem]{Example}

\numberwithin{theorem}{section}

\newcommand{\Z}{{\mathbb{Z}}}
\newcommand{\Q}{{\mathbb{Q}}}
\newcommand{\F}{{\mathbb{F}}}

\newcommand{\R}{{\mathbb{R}}}

\newcommand{\C}{{\mathbb{C}}}
\newcommand{\trace}{\mbox{trace}}

\renewcommand{\em}{\sf}

\begin{document}
 \bibliographystyle{plain}

\title{Orthogonal determinants of  characters}
\author{Gabriele Nebe}
\email{nebe@math.rwth-aachen.de}
\address{Lehrstuhl f\"ur Algebra und Zahlentheorie, RWTH Aachen University, 52056 Aachen, Germany}

\begin{abstract}
 For an irreducible orthogonal character $\chi $ of even degree 
	there is a unique square class $\det({\chi })$ in the character field 
 such that the invariant quadratic forms in any $L$-representation 
	affording $\chi $ have determinant in $\det({\chi })(L^{\times })^2$. 
	\\
	{\sc keywords:} orthogonal representations of finite groups; Schur index.
\end{abstract}

\maketitle

\section{Introduction}
\label{sec:intro}

 Let $G$ be a finite group. 
 An absolutely irreducible complex representation 
 $\rho _{\C }:G \to \GL_n(\C) $ fixes a non-zero quadratic form
 if and only if $\rho _{\C } $ is equivalent to a real representation
 $\rho : G\to \GL_n(\R) $. 
 Then $B:= \sum _{g\in G} \rho(g) \rho(g)^{tr} $ is a 
 positive definite symmetric matrix. As $\rho (g) B \rho(g)^{tr}  = B$
 the group $\rho (G)$ is a subgroup of the orthogonal group of $B$.
 In this case we call $\chi : G\to \R , 
 \chi (g) := \trace (\rho (g)) $  
 an {\em orthogonal character}. 
 Let $K= \Q (\chi (g) \mid g\in G) $ be the {\em character field} of $\chi $.
 The main result of this paper is Theorem \ref{main} 
 showing that for even character degree 
 there is a unique square class $\det({\chi } )\in K^{\times }/(K^{\times })^2$
 such that for any representation $\rho $ with character $\chi $ 
 over some extension field $L$ of $K$ all non-zero $\rho(G)$-invariant
 symmetric bilinear forms $B$ have $\det(B) \in \det({\chi }) (L^{\times })^2$. 
 The square class $\det({\chi })$ is called the {\em orthogonal determinant} 
 of $\chi $.
The proof is immediate when the Schur index of $\chi $ 
is one. In this case there is a representation $\rho $ for $L=K$ and
 $\det({\chi }) = \det(B) (K^{\times })^2$ for any non-zero $\rho(G)$-invariant form $B$. 
If the Schur index of $\chi $ is two, 
there is no such representation over $K$. 
Then the rational span of the matrices in $\rho (G)$ is 
a central simple $K$-algebra $A$ and 
by Remark \ref{adjointinverts} 
the adjoint involution induces an involution on $A$ 
whose determinant (as defined in \cite[Proposition (7.1)]{Tignol}) 
is $\det({\chi })$.

In positive characteristic all Schur indices are one 
 and the result of Theorem \ref{main} holds with a direct  
 easy proof. Therefore we restrict to characters over number fields 
 in this short note. 


 In an ongoing project with Richard Parker, 
 we aim to provide the orthogonal determinant 
 for all irreducible orthogonal Brauer characters 
for all but the largest few ATLAS groups \cite{ATLAS} over all finite fields. 
 This is a finite (computational) problem for 
 the primes that divide the group order.
 Thanks to Theorem \ref{main}
 the infinitely many primes not dividing the group order 
 can be treated with 
 a characteristic zero approach as 
 illustrated in Corollary \ref{J2}. 

 {\sc Acknowledgements.}
I thank Richard Parker for his persisting questions and 
fruitful discussions that made me write this paper and 
Eva Bayer-Fluckiger for her interesting comments leading to
Section \ref{appl}.


%
%

\section{Determinants of symmetric bilinear forms} 

Let $K$ be a field of characteristic $\neq 2$,
$V$ a vector space  of dimension $n$ over $K$ 
and $\tilde{B}: V\times V \to K$ a symmetric bilinear form.
Any choice of a basis
 $(e_1,\ldots , e_n)$ of $V$ identifies $V$ with the row space $K^{ n}$.
The Gram matrix of $\tilde{B}$ with respect to this basis is
$$B:= (\tilde{B}(e_i,e_j)) _{i,j=1}^n \in K^{n\times n}$$ 
a symmetric square matrix satisfying 
$ \tilde{B} (x,y) = x B y^{tr} $
for all $x,y\in K^{n}$. 
Base change by the matrix $T \in \GL_n(K)$ changes the Gram matrices into 
$T B T^{tr}$ and hence 
the determinant of $\tilde{B}$ is 
 $$\det(\tilde{B}) :=\det (B) (K^{\times })^2 \in K/(K^{\times })^2$$ 
 well defined up to squares. 
 The
  bilinear form $\tilde{B}$
 is called {\em non-degenerate}, if $\det(\tilde{B}) \in K^{\times }/(K^{\times })^2$,
i.e. $\det (B) \neq 0$. 

\subsection{The adjoint involution}

Any non-degenerate symmetric bilinear form $\tilde{B}$ on $V$
defines a $K$-linear involution $\iota_{\tilde{B}}$ on $\End_K(V)$. 
For $\alpha \in \End_K(V)$ the endomorphism $\iota_{\tilde{B}}(\alpha )$ is defined by
$$ \tilde{B}(\alpha(x),y) = \tilde{B}(x,\iota _{\tilde{B}}(\alpha) (y)) \mbox{ for all } x,y\in V.$$
Identifying $\End_K(V)$ with the matrix ring $K^{n\times n}$ 
by choosing a basis $(e_1,\ldots, e_n)$ of $V$, the involution 
$\iota _{\tilde{B}}$ is given by 
$$\iota _B(A) = B A^{tr} B^{-1} .$$ 
We define 
$$E_-({\tilde{B}}) := \{ \alpha \in \End_K(V) \mid \iota _{\tilde{B}} (\alpha ) = - \alpha \} $$ 
the $K$-space of skew adjoint endomorphisms.
In matrix notation we get
$$E_-(B):= \{ X\in K^{n\times n} \mid B X^{tr} B^{-1} = -X \} $$ 
and hence 

\begin{lemma} \label{a}
$E_-(B) = \{ BX \mid X = -X^{tr} \in K^{n\times n } \} $.
\end{lemma} 


Scaling of the bilinear form does not change the involution,
$E_{-} (aB) = E_{-} (B)$ for all $a\in K^{\times }$. 
On the other hand $\det (aB) = a^n \det(B)$.
So we can only read off the determinant of $\tilde{B}$ from the involution 
$\iota _B$ in even dimensions. 
The following property of skew adjoint endomorphisms is crucial.

\begin{proposition} \label{detB}
	$E_-(\tilde{B})$ contains  invertible elements if and only if $\dim(V)$ is even. 
	Then $\det(\tilde{B}) =  \det (\alpha ) (K^{\times })^2 $ for 
any invertible $\alpha \in E_-(\tilde{B})$.
\end{proposition}

\begin{proof}
	We prove the theorem in matrix notation. 
	Let 
$E_-(I) := \{ X\in K^{n\times n } \mid X = -X^{tr} \} $ denote the
	space of skew symmetric matrices.
It is well known that $E_-(I)$ contains an invertible 
matrix, if and only if $n$ is even and then the determinant
of such a matrix is a square. 
By Lemma \ref{a}  the 
map $E_-(I) \to E_-(B) , X\mapsto  B X$ is an isomorphism. 
So $E_-(B) $ contains invertible elements 
if and only if $\dim(V)$ is even, and all 
	such elements $Y\in E_-(B) \cap \GL_n(K)$ satisfy $\det(Y) \in \det(B)(K^{\times })^2$.
\end{proof} 

\subsection{Determinants and isometries} \label{appl} 

For any non-degenerate symmetric bilinear form $\tilde{B}$ its 
{\em orthogonal group} is
$$O(V,\tilde{B}):= \{ g\in \GL(V) \mid \tilde{B}(vg,wg)=\tilde{B}(v,w) \mbox{ for all } v,w\in V \} .$$

\begin{remark}
	An endomorphism $g\in \End _K(V)$ lies in $O(V,\tilde{B})$ if and only if 
	$g \iota _{\tilde{B}}(g) = \iota_{\tilde{B}}(g) g = \id _V$, i.e. $\iota_{\tilde{B}}(g) =g^{-1}$. 
\end{remark}

\begin{proposition} (see \cite[Proposition 5.1]{Eva}) \label{cyc}
	Let $g\in O(V,\tilde{B})$ and denote by $P$ the characteristic polynomial
of $g$. Assume that $P(1)P(-1) \neq 0$.
\begin{itemize} 
\item[(a)] $\dim (V) $ is even.
\item[(b)] $\det (g) = 1$.
\item[(c)] $\det(\tilde{B}) = \det(g-g^{-1}) (K^{\times})^2 = P(1) P(-1) (K^{\times })^2 $.
\end{itemize}
\end{proposition} 

\begin{proof}
The (sketched) proof 
follows the exposition in \cite{Eva}. 
 \\
	Let $n:=\dim (V) = \deg(P)$. Then $P(0) = (-1)^n \det(g) =: \epsilon \in \{ 1, -1 \} $. 
	\\
	Put
	$P^*(X) := \epsilon X^n P(X^{-1})$ to denote the reverse polynomial
	of $G$. By \cite[Proposition 1.1]{Eva} the condition that
	$g\in O(V,\tilde{B})$ implies that $P=P^{*}$. 
	As $P(1) \neq 0$ we hence have $\epsilon = 1$. Now $P(-1)\neq 0$ 
	yields that $n$ is even and so $\det(g) = 1$.
	\\
	To see (c) we write $P(X) = \prod_{j=1}^n(X-\xi_j) $ over 
	some algebraic closure of $K$. 
	Then 
	$$P(1)P(-1) =  \prod_{j=1}^n(\xi_j^2 -1) = (\prod_{j=1}^n \xi_j) 
	\prod_{j=1}^n (\xi_j-\xi_j^{-1}) = \det(g) \det(g-g^{-1}) .$$
	As $\det(g) =1$ and $g-g^{-1} \in E_-(\tilde{B}) $ is a unit,
	statement (c) now follows from Proposition \ref{detB}.
\end{proof}

\begin{corollary} \label{det1} 
	If there is $g\in O(V,{\tilde{B}})$ with $g^2 = -1$ then $\det({\tilde{B}}) = 1$.
\end{corollary}

\begin{proof}
The minimal polynomial of $g$ divides $X^2+1$
and hence the characteristic polynomial of $g$ is $P=(X-i)^a(X+i)^b$ with
 $P(1)P(-1) = (-2)^{a+b}$.
 By Proposition \ref{cyc} (a) the degree of $P$ is $a+b = \dim (V) $ is even so 
 $P(1)P(-1)$ is a square and 
	the statement follows from Proposition \ref{cyc} (c).
\end{proof}

\section{Orthogonal representations of finite groups}

Let $G$ be a finite group and $L$ be a field. 
An $L$-representation $\rho $ is a 
group homomorphism 
$\rho : G \to \GL_n(L)$.  Given a representation $\rho $ we put
$${\mathcal F}(\rho ) := \{ B\in L^{n \times n } 
\mid  B = B^{tr} \mbox{ and } 
\rho(g) B \rho (g)^{tr} = B \mbox{ for all } g\in G\} $$ to denote
the $L$-vector space of symmetric $G$-invariant bilinear forms on $L^{n}$. 

\begin{remark}\label{adjointinverts}
	Let 
	$B\in {\mathcal F}(\rho ) \cap \GL_n(L)$. 
Then the adjoint involution 
$\iota_B$ on $L^{n\times n}$ satisfies $\iota _B(\rho(g)) = \rho (g^{-1} )$ 
for all $g\in G$.
\end{remark}

\begin{definition}
	A representation $\rho $ and also its character $\chi_{\rho } = \trace \circ \rho $ is called {\em orthogonal}, if ${\mathcal F}(\rho )$ contains 
	a non-degenerate element. 
\end{definition}

\subsection{Orthogonal determinants}

\begin{theorem}\label{main}
Let  $\chi \in \Irr_{\C} (G)$ be  an orthogonal character 
 of even degree $n:=\chi (1) \in 2\Z $. Denote by $K$ the 
character field of $\chi $. 
Then there is a unique totally positive 
	square class $\det({\chi }) = d (K^{\times })^2 \in K^{\times }/(K^{\times })^2 $ 
with the following property:
Let $L\supseteq K$ be any field 
and
$\rho  : G\to  \GL _{n} (L) $ 
 a representation 
with character $\chi $.  Then  any nonzero 
$B\in {\mathcal F}(\rho )$ satisfies
	$\det (B) \in \det({\chi })(L^{\times })^2$. 
\end{theorem} 

\begin{definition}
	The square class $\det({\chi }) (K^{\times })^2 $ from Theorem \ref{main} 
	is called the {\em orthogonal determinant} of the character $\chi $.
\end{definition}

%

\subsection{Proof of Theorem \ref{main}}\label{proof}

For the proof of Theorem \ref{main} we assume that we are 
given a field $L$ containing the character field $K$ of $\chi $ and 
a representation $\rho : G\to \GL_n(L)$ with character $\chi $. 
We also choose some non-zero $B \in {\mathcal F}(\rho )$.
Since $\rho $ is absolutely irreducible 
the matrices in $\rho(G)$ generate $L^{n\times n}$ as 
a vector space over $L$. Also ${\mathcal F}(\rho )$ 
is one dimensional and $B$ is non-degenerate and unique up to scalars:

\begin{remark}
	\begin{itemize}
		\item[(a)]
	$\langle \rho (g) \mid g \in G \rangle _L = L^{n\times n}$
\item[(b)]
 ${\mathcal F}(\rho ) = \langle B \rangle _{L}$.
 \item[(c)] 
	 $E_-(B) = \langle \rho(g) - \rho (g^{-1}) \mid g\in G \rangle _{L}$.
	\end{itemize}
\end{remark}

We now consider the $\Q $-algebra generated by the matrices in $\rho (G)$,
$$A := \langle \rho (g) \mid g\in G \rangle _{\Q } \leq L^{n\times n}.$$
\begin{remark} 
	\begin{itemize}
		\item[(a)]
$A$ is a central simple $K$-algebra of dimension $n^2$. 
\item[(b)]
	The restriction $\iota $ of $\iota _B$ to $A$ 
			satisfies $\iota(\rho (g)) = \rho(g^{-1})$ for all $g\in G$.
		\item[(c)]
			$E_-(\rho ) := \langle \rho(g) - \rho(g^{-1}) \mid g\in G \rangle _{\Q} =  E_-(B) \cap A$.
	\end{itemize}
\end{remark}

It is well known that the reduced norm of a central simple algebra 
takes values in the center of this algebra:

\begin{lemma} \label{Nred}
For all $X\in A$ we have that $\det(X) \in K$.
\end{lemma} 

\begin{proof}
	As $L$ is a splitting field for $A$ the determinant is the 
	reduced norm of the central simple $K$-algebra $A$, see for instance
	\cite[Section 9]{Reiner}. 
	Reiner also shows that 
the reduced norm is independent of 
	the choice of a splitting field and takes values in $K$. 
\end{proof}

	By \cite[Corollary 2.8]{Tignol} a central simple algebra 
	with orthogonal involution contains invertible elements that are 
	negated by the involution if and only if the dimension of this
	algebra over its center is even.  
	In particular for our situation this yields the following 
	proposition
 for which we give an independent short proof below.

\begin{proposition}  \label{skewinA}
	$E_-(\rho )$ contains  invertible elements.
\end{proposition}

\begin{proof}
	The fact that $E_-(\rho )$ contains an element that 
	is invertible in the central simple $K$-algebra 
	$A = \langle \rho (g) \mid g\in G \rangle _{\Q }$ 
	does not depend on the choice of 
	the splitting field $L$.
	So without loss of generality we fix an embedding 
	$\epsilon : K \hookrightarrow \R $, identify $K$ with its image $\epsilon (K) \subseteq \R $, and take $L=\R $, one of the real completions of $K$.

	We first choose a $K$-basis $B=(b_1,\ldots , b_m)$ of 
	$E_-(\rho )$. Then $B$ is also an $\R$-basis of $E_-(B)$. 
	Let $X\in E_-(B ) \cap \GL_n(\R)$ be an invertible element 
	of $E_-(B)$. 
	Write 
	$$ X = \alpha _1 b_1+ \ldots + \alpha _m b_m \mbox{ with unique } \alpha _i \in \R .$$
	Now $\Q $ and hence also $\epsilon(K)$ is dense in $\R $.
	So there are $a_i \in K $ such that $\epsilon (a_i) $ 
	is arbitrary close to $\alpha _i$ for all $i=1,\ldots , m$. 
	Put
	$$ Y := a_1 b_1 +\ldots + a_m b_m \in E_-(\rho )  .$$
	For $Y$ being a unit in $A$,  it is enough to achieve that  the determinant of 
	$\epsilon(Y) := \sum _{i=1}^m \epsilon (a_i) b_i $ is non zero. 
	As $\det $ is a polynomial, in particularly continuous, 
	and $\det(X) \neq 0$, we can find $a_i \in K$ 
	such that $\epsilon(\det(Y)) = \det(\epsilon (Y)) \neq 0$. 
	But then $Y \in E_-(\rho ) $ is an invertible matrix.
\end{proof}

\begin{proof}
	(of Theorem \ref{main}) By Proposition \ref{detB} 
	we get $\det (B) (L^{\times })^2 = \det (X) (L^{\times })^2$ for any invertible 
$X\in E_-(B)$. 
Proposition \ref{skewinA} says that such an invertible
element $X$ can be chosen in 
$E_-(\rho ) = E_-(B)\cap A$, so in particular its determinant 
is an element of $K$ by Lemma \ref{Nred}. 
\end{proof}

\section{Some applications}

\subsection{An example: $\SL_2(\F_7)$}

For illustration let $G:=\SL_2(\F_7)$ be 
the special linear group of degree 2 over the field with 7 elements. 
The complex character table of $G$ is given in \cite{ATLAS}.
For any faithful irreducible representation $\rho $ of 
$G$ we obtain $\rho (\left(\begin{array}{cc} 0 & -1 \\ 1 & 0 \end{array} \right) )^2 = -\id $ and hence 
by Corollary \ref{det1} all 
 faithful irreducible orthogonal characters have determinant 1.
There are 
 six complex irreducible characters of the group $L_2(7) = \PSL_2(7)$,
 of degrees 
$1$, $3,3$, $6$, $7$, and $8$ 
giving rise to three irreducible rational representations of 
even degree, $3ab$, $6$, and $8$:
	\begin{itemize} 
		\item[3ab]
			Restrict the representation 
			to a Sylow-7-subgroup $\langle g\rangle$ of $G$. 
			The eigenvalues of $\rho(g) $ are all 
			primitive 7th roots of unity and hence 
			$\det(\rho(g) - \rho(g^{-1}) ) =
			\prod_{i=1}^6 (\zeta _7^i - \zeta _7^{-i}) = 7$.
			So by Proposition \ref{cyc} the orthogonal 
			determinant of $\chi $ is $\det({\chi }) = 7 (\Q^{\times })^2$.
		\item[6] Restriction to  
			$\langle g \rangle $ as before allows to
			conclude that $\det({\chi} )=7 (\Q^{\times })^2$.
		\item[8] 
Let $H = C_7:C_3 = \langle g \rangle : \langle h \rangle$ be the 
normaliser in $G$ of the Sylow-7-subgroup.
Then the restriction of $\rho $ to $H$ decomposes as $6+2$, 
where $\langle g \rangle $ acts fixed point free on the 
$6$-dimensional part and trivially on the $2$-dimensional summand. 
On the $2$-dimensional summand, $\langle h \rangle $ acts faithfully.
If $e = \frac{1}{7} \sum _{i=0}^6 g^i  \in \Q H$ is the involution 
invariant idempotent projecting onto the fixed space of $\langle g \rangle$, 
then 
$$X:= (\rho(g) -\rho(g^{-1})) + \rho(e) (\rho (h) - \rho (h^{-1}) )  \in E_-(\rho ) $$ 
			has determinant $7\cdot 3 = 21$. So by Proposition \ref{detB} we get $\det({\chi }) = 21 (K(\chi )^{\times })^2$.
	\end{itemize} 

\begin{remark}
Note that these techniques essentially suffice to find all 
orthogonal determinants 
for all groups $\SL_2(q)$ as given in \cite{Braun}.
\end{remark}

\subsection{Orthogonal characters with rational Schur index 2}

Theorem \ref{main} is particularly helpful in the case that 
the orthogonal character  is not the character of 
a representation over its character field. 
The smallest example of a simple group $G$ in \cite{ATLAS} is 
the group $G=J_2$. 
This sporadic simple group 
has a complex irreducible orthogonal character $\chi $ of degree 
$\chi (1) = 336$ with rational character field.
By \cite{Feit} (see also \cite{Unger}) the rational Schur index of $\chi $ 
is $2$.
With MAGMA \cite{MAGMA} 
we realise the representation as a rational representation $\rho $
of dimension $2\cdot 336 = 672$ (with character $2\chi $). 
Then 
the central simple $\Q $-algebra $$A=\langle \rho (g) \mid g\in G \rangle 
\cong {\mathcal Q}_{2,3} ^{168 \times 168} $$
is isomorphic to a matrix ring over the indefinite rational quaternion algebra 
${\mathcal Q}_{2,3} $ ramified at 
$2$ and $3$. 
We take three random elements $g_1,g_2,g_3\in G$ to achieve 
that $x:=\sum_{i=1}^3 \rho(g_i) - \rho(g_i^{-1}) \in E_-(\rho )$ 
has full rank. 
To compute the reduced norm of $x \in A$, we compute 
the characteristic polynomial $P$ of $x$ which is a square $P=p^2$ 
of a unique monic polynomial $p$. 
Then the reduced norm of $x$ is $p(0)$. 
It turns out that $p(0)$ is a rational square,
so $\det({\chi }) = 1$.

\begin{corollary} \label{J2}
	For any finite field $F$ of characteristic $p \geq 7$ 
	the representation $\rho : J_2 \to \GL_{336}(F) $ with 
	Brauer character $\chi $ fixes a symmetric bilinear form of determinant 1.
	In particular $\rho (J_2) \leq O^+_{336}(F)$.
\end{corollary}

\subsection{Split extensions $G:2$} 

Let $G$ be a finite group, $\alpha \in \Aut (G)$ an automorphism of order 2.
Then the split extension $G:2$ has a pseudo-presentation 
$$G:2 = \langle G,h \mid hgh^{-1} = \alpha (g), h^2 =1 \rangle .$$
Assume that there is an orthogonal  character $\chi \in \Irr_{\C }(G)$ 
such that $\chi \circ \alpha \neq \chi $. 
Then there is a unique character $\Chi \in \Irr _{\C}(G:2)$ 
such that $\Chi _{|G} = \chi + \chi \circ \alpha $. 
As $\Chi (hg) = 0$ for all $g\in G$ the 
character field $F$ of $\Chi $ is contained in the character field 
$K$ of $\chi $. 

\begin{theorem}\label{induz}
	If $K=F$ then $\det (\Chi ) = 1$. \\
	Otherwise $K$ is a quadratic extension of $F$, so $K=F[\sqrt{\delta }]$ 
	for some $\delta \in F$ and $\det(\Chi ) = \delta^{\chi(1)}$.
\end{theorem}

\begin{proof}
	Let $L$ be some extension of $K$ and $\rho : G\to \GL_n(L)$ 
	a representation affording the character $\chi $, so $n=\chi (1)$.
	Then the induced representation $\Rho $ with 
	character $\Chi $ is given by 
	$$\Rho(g) =  \diag (\rho (g) , \rho (\alpha (g)) ) \mbox{ for all }  g\in G \mbox{ and } \Rho (h) =  
\left( \begin{array}{cc} 0 & 1 \\ 1 & 0 \end{array} \right) .$$
	In particular $\Rho $ is also an $L$-representation and 
	$${\mathcal F}(\Rho ) = \{ \diag (B,B) \mid  B\in {\mathcal F}(\rho ) \} .$$
	This shows that $\det(\Chi ) (K^{\times })^2 = \det(\chi )^2 $.
	In particular $\det(\Chi )  = 1$ if $K=F$. 
	\\
	Now assume that $K=F[\sqrt{\delta}]$. 
	Let $g_0\in G$ be such that
$K=F[\chi(g_0) ]$ and put $C_0:= \sum _{g\in g_0^G} g $ the
class sum of $g_0$. Then
$C_0$ and $\alpha (C_0)$ are central elements in $\Q G$.
Adding some element of $F$ and
	multiplying by some element in $F^{\times }$ we
find $C$ in the center of $\Q G$, such that
	$$ \chi (C)= \chi (\iota(C)) = n \sqrt{\delta } , \chi (\alpha (C)) = \chi (\alpha (\iota(C))) = -n \sqrt{\delta } .$$
	As $h=h^{-1}$ and also $\iota (C) = C $ we compute
	$$\begin{array}{lll} \iota(\Rho(h) \Rho(C)) = & \Rho(\iota(C)) \Rho(h) =
		& 
		\Rho (C) \Rho (h) = \\  \Rho(h) \Rho(\alpha (C)) = &  \Rho(h) (-\Rho(C)) = &
	-\Rho(h) \Rho(C) .\end{array} $$
		    So $\Rho(h) \Rho(C) \in E_-(\Rho )$ and
 $\det(\Rho(h)\Rho(C)) = (-1)^n \sqrt{\delta }^n (-\sqrt{\delta } ) ^n = \delta ^n $. 
	Thanks to Proposition \ref{detB} and 
	Theorem \ref{main} we get  $\det({\Chi }) = \delta ^n (F^{\times })^2$.
\end{proof}


\begin{example}
	Let $\Chi _{n} \in \Irr _{\C } (J_2:2) $ be the irreducible 
	characters of degree $2n$  of the automorphism group of $J_2$, 
	for $n=14,21,70,189,224 $ (see \cite{ATLAS}). 
	In all cases the character field of $\Chi _n$ is $\Q $ and 
	the restriction of $\Chi _n$ to the simple group $J_2$ is the 
	sum of two irreducible orthogonal characters of degree $n$ and
	with character field $\Q [\sqrt{5}]$. 
As $J_2:2$ is a split extension  Theorem \ref{induz} tells us that 
	$\det(\Chi _n)= 5^n (\Q ^{\times })^2$. 
\end{example}

\end{document}